\documentclass[11pt, oneside]{article}
\usepackage{amsmath,amsfonts,amssymb,amscd}
\usepackage{amsthm}

\usepackage{eucal}

\date{}

\begin{document}

\bibliographystyle{plain}

\newtheorem{theorem}{Theorem}[section]
\newtheorem{proposition}{Proposition}[section]
\newtheorem{lemma}{Lemma}[section]
\newtheorem{corollary}[lemma]{Corollary}

\title{Manifolds of positive Ricci curvature with quadratically asymptotically nonnegative curvature and infinite topological type}

\author{Huihong Jiang and Yi-Hu Yang\thanks{Partially supported by NSF of China (No.11571228)}}

\maketitle

\begin{abstract}

We construct a complete $n$-dimensinal ($n \ge 6$) Riemannian manifold of positive Ricci curvature with
quadratically asymptotically nonnegative sectional curvature and infinite topological type. This gives a negative answer
to a problem proposed by Jiping Sha and Zhongmin Shen \cite{SS} in the case of $n \ge 6$.
\end{abstract}

\section{Introduction}

The purpose of the paper is to construct a noncompact complete Riemannian manifold of nonnegative Ricci curvature
with quadratically asymptotically nonnegative sectional curvature and infinite topological type;
this gives a negative answer to the problem proposed by Jiping Sha and Zhongmin Shen \cite{SS},
which asks {\it if a complete Riemannian manifold with nonnegative Ricci curvature
and quadratically asymptotically nonnegative sectional curvature is of finite topological type}.

In order to give a more precise statement, let's first fix some notation. Let $M^n$ be an $n$-dimensional
complete noncompact manifold, $p_0\in M$ a fixed point. By $B(p_0, t)$ (resp. $S(p_0, t)$)
denote the ball (resp. sphere) with the center at $p_0$ and the radius $t$. Set
$$
{\text{diam}}(p_0; t)=\sup_i{\text{diam}}(\Sigma_i, M\setminus B(p_0, {\frac 12}t)),
$$
where $\Sigma_i$ is a connected component of $S(p_0, t)$, ${\text{diam}}(\Sigma_i, M\setminus B(p_0, {\frac 12}t))$
denotes the interior diameter of $\Sigma_i$ in $M\setminus B(p_0, {\frac 12}t)$ when $M\setminus B(p_0, {\frac 12}t)$
is considered as a metric space with the induced metric and $\Sigma_i$ as its subset.
We also define
$$
K_{p_{0}} (t) = \inf_{M^{n} \setminus B(p_{0},t)} K,
$$
where $K$ denotes the sectional curvature of $M$, and the infimum is taken over all the 2-planes at the points
in $M\setminus B(p_{0},t)$. We say that $M$ is of {\it quadratically asymptotically nonnegative curvature} if,
for all $t\geq 0$,
$$
K_{p_{0}} (t) \geq -{\frac{K_{0}}{1+t^2}},
$$
for some positive constant $K_0$.
We also say that $M$ is of finite topological type if it is homeomorphic to
the interior of a compact manifold with boundary; otherwise, infinite topological type.
Then, the main aim of this paper is to show the following

\begin{theorem} There exists a $6$-dimensional complete Riemnannian manifold $M$
with a base point $p_0$, which is of positive Ricci curvature, quadratically asymptotically nonnegative curvature,
and infinite topological type. Furthermore, it satisfies
\begin{align*}
&\lim_{t\to\infty}{\frac{{\text{vol}}(B_{t}(p_0))}{t^6}}=0~~{\text{and}}~~{\text{diam}}(p_{0};t)=O(t).
\end{align*}
\end{theorem}

Earlier examples of positive (or nonnegative) Ricci curvature and infinite
topological type have been constructed by Sha and Yang \cite{ShY},
which have sectional curvature bounded below but diameter growth of order
$O(t^{\frac 2 3})$ (also see \cite{SW}); moreover, they have infinite second Betti number.
The examples of Sha and Yang \cite{ShY, ShY1} show that in the higher dimensional case ($\ge 4$)
one cannot expect similar results of Cheeger-Gromoll \cite{CG}
in the case of positive (or nonnegative) Ricci curvature.
In other words, in order to get certain finiteness result,
some additional assumptions (e.g. volume growth or diameter growth) must be required;
in this aspect, one can refer to
\cite{AG, SW, SS, LS, dx, sso, JY} and therein. As for the 3-dim case, complete noncompact Riemannian
manifolds with nonnegative Ricci curvature were completely classified recently by G. Liu \cite{L} based on
the technique earlier developed by Schoen-Yau \cite{SY}.

After the examples of Sha and Yang, Perelman \cite{P} gave a more general construction
of 4-dimensional compact examples with positive Ricci curvature, big volume and arbitrary
large Betti numbers. Menguy applied Perelman's construction to obtain noncompact examples
of infinite topological type \cite{M1,M2}. However, all these noncompact examples are
not of quadratically asymptotically nonnegative curvature though they have positive
Ricci curvature and infinite topological type. In particular, Menguy's examples \cite{M1} have
sectional curvature bounded below by $-({\frac{g(t)} {t}})^2$; here, $t$ is the distance function
to a fixed point and $g(t)$ is a function that is going to infinity as $t$ goes to infinity.

It is easy to get higher dimensional examples with the same properties by taking
a metric product of the above 6-dimensional example with the standard spheres.
Here, we should also point out that the examples with only quadratically
asymptotically nonnegative curvature and infinite topological type has been constructed
by Abresch \cite{A}; but these examples are not of nonnegative Ricci curvature in general.

In the final of this section, we give a brief description for the idea of the construction,
which mainly benefits from the works of Perelman \cite{P} and Menguy \cite{M1, M2, M3}.
We first begin with a doubly warped product metric
$$
ds^{2}=dt^{2}+u^{2}(t)(dx^{2}+f^{2}(t,x)d\sigma^{2})+g^{2}(t)d\theta^{2},
$$
where $d\sigma^{2}, d\theta^{2}$ are the standard metric of the 2-dim sphere $S^{2}$.
For the first part, $dt^{2}+u^{2}(t)(dx^{2}+f^{2}(t,x)d\sigma^{2})$, of the metric,
we assume that it is the spherical metric for $t_{i} < t < t_{i}+2r_{i}$, here $t_i+2r_i<t_{i+1}$
(for the choice of $t_i$ and $r_i$, cf. \S 2), and $u(t)=O(t)$.
Thus, the metric on this part is similar to the compact example of Perelman \cite{P},
which is a double suspension over a small 2-sphere so that one gets a singular 4-sphere
with the set of its singular points being a circle. Then, it allows us to remove the geodesic ball
$B_{\frac 45 r_{i}}(o_{i})$($o_{i}=(t_{i}+r_{i},0)$) and glue in a $\mathbb{CP}^{2}$
via a neck, provided with $f(t,x)$ almost like $R_{0}\sin x$ ($R_{0} \ll 1$)
and smoothed near $x=0$. We can continue this surgery when $i$ goes to infinity.
Then, we get a manifold with infinite second Betti number, which is of course infinite topological type.

Actually, this is exactly the constructing process of Menguy's Example with Euclidean
volume growth \cite{M1} (i.e. the order of growth is the dimension of manifolds).
However, the radius $r_{i}$ in Menguy's examples must be of order $o(t)$
(more precisely, $r_{i}=\frac{t_i}{h(t_i)}$, here $h$ is a certain function going to infinity
as $t\to\infty$), this makes the sectional curvature to
be bounded from below by $-({\frac{h(t)}{t}})^2$ for sufficiently large $t$, which certainly is not
quadratically asymptotically nonnegative. But the factor $h(t)$ here 
is essential for $-\frac{u_{tt}}{u} \ge 0$ so that Ricci curvature can be
controlled to be positive. In our construction, we set $r_{i}=rt_i=O(t)$ to get
quadratically asymptotically nonnegative sectional curvature (here $r$ is a certain positive constant, see \S 2.1).
To ensure that Ricci curvature is positive, we take the second warped product part $\times_{g(t)}S^2$
and assume the order $\gamma$ of $g(t)$ ($=O(t^{\gamma})$) being sufficiently small,
so that we can control the curvature through some new directions.
Moreover, in order to be able to do surgery by gluing $\mathbb{CP}^{2}$ via a neck on
$(t_{i}+\frac{r_{i}}{6}, t_{i}+\frac{11r_{i}}{6})$,
we set $g(t)$ to be constant on $(t_{i}+\frac{r_{i}}{6}, t_{i}+\frac{11r_{i}}{6})$.
We here remark that since the first warped product on $[t_i, t_i+2r_i]$ is part of a certain sphere
(i.e. $u(t)={\frac{1}{\sqrt{K_i}}}\sin(\sqrt{K_i}(t-t_i+\psi_i))$, see \S 2.1), our choice for $r_i$
will not destroy the positiveness of Ricci curvature on this part;
so, we can assume that $g$ is constant on $(t_{i}+\frac{r_{i}}{6}, t_{i}+\frac{11r_{i}}{6})$.
Clearly, by means of Sha-Shen's finiteness result \cite{SS}, the manifold constructed here
is not of Euclidean volume growth.

\section{Construction of the manifolds}

We begin with a doubly warped product manifold with boundary,
$$
Q = [t_1,+\infty) \times_{u} S^3 \times_{g} S^2,
$$
for some $t_1>0$ and certain functions $u(t)$, $g(t)$ in $t\ge t_1$, all of which will be fixed later.
Note that $S^3$ here is not necessarily a round sphere.

First, taking a sequence of intervals $[t_i, t_i +2 r_i]$ ($1 \leq i < +\infty$) to be determined later,
we will let $u(t)$ be spherical on $[t_i, t_i +2 r_i]$ and $g(t)$ be constant
$g_i$ on $[t_{i}+\frac{r_{i}}{6}, t_{i}+\frac{11r_{i}}{6}]$.
In the following, we'll see that our surgeries actually arise on
$(t_{i}+\frac{r_{i}}{6}, t_{i}+\frac{11r_{i}}{6})\times_{u} S^3 \times_{g_i} S^2$, $1 \leq i < +\infty$.

In fact, we have surgeries only on $(t_{i}+\frac{r_{i}}{6}, t_{i}+\frac{11r_{i}}{6})\times_{u(t)} S^3$.
More precisely, we remove a geodesic ball $B_{\frac 45 r_{i}}(o_{i})$
with the center at $o_{i}=(t_{i}+r_{i},0)$ (for the precise meaning of $o_i$, see \S 2.1)
on $(t_{i}+\frac{r_{i}}{6}, t_{i}+\frac{11r_{i}}{6}) \times_{u(t)} S^3$,
and glue a $\mathbb{CP}^{2} \setminus B(p_0)$ (called a {\it core}, topologically $\mathbb{CP}^2$ deleting a ball)
via $S^3 \times [0,1]$ (called a {\it neck}) along $\partial B_{\frac 45 r_{i}}(o_{i})$ with appropriate metrics depending on $i$
and the boundary properties on $\mathbb{CP}^{2} \setminus B(p_0)$ and $S^3 \times [0,1]$ (for details, see \S 2.2).
By $P_i$ ($1 \leq i < +\infty$) denote these glued parts (both the core and the neck).
Similarly, we also glue a core via a neck along the boundary $t=t_1$ to smooth the manifold near the origin.
We denote this glued part by $P_0$.

So, our manifold is actually
$$
M = (Q \setminus \coprod\limits_{i=1}^{+\infty} (B_{\frac 45 r_{i}}(o_{i}) \times_{g_i} S^2))
\cup_{\textrm{Id}} \coprod\limits_{i=0}^{+\infty} (P_i \times_{g_i} S^2)
$$
where $g_0 = g_1$ and 'Id' means gluing along the boundary through the identity map.

Here, we need to point out that actually one should also do the same surgery
on $Q$ at $o'_{i}=(t_{i}+r_{i},\pi)$, i.e. remove a geodesic ball $B_{\frac 45 r_{i}}(o'_{i})$
(isometric to $B_{\frac 45 r_{i}}(o_{i})$) with the center at $o'_{i}=(t_{i}+r_{i},\pi)$
and glue with the $P_i$ along the boundary, since the Ricci curvature is also not necessarily positive
on $B_{\frac 45 r_{i}}(o'_{i})$. In fact, our construction of $Q$ in $x$ is symmetric at $x={\frac \pi 2}$. More precisely,
$f(t,x)$ will first be constructed for $x\in [0,\frac\pi2]$ and then symmetrically extended to $[\frac\pi2,\pi]$;
moreover $B_{\frac 45 r_{i}}(o_{i})$ and $B_{\frac 45 r_{i}}(o'_{i})$ are disjoint (for details, see \S 2.1.3).
Based on this, for the sake of simplicity, in the following we will only do the surgery on $B_{\frac 45 r_{i}}(o_{i})$.

$M$ has infinite second Betti number and hence be of infinite topological type.
In the following, we will give much more details of the construction.

\subsection{Construction of $Q$}
We equip $Q=[t_1,+\infty) \times_{u(t)} S^3 \times_{g(t)} S^2$ with the following metric
$$
ds^{2}=dt^{2}+u^{2}(t)[dx^{2}+f^{2}(t,x)d\sigma^{2}]+g^{2}(t)d\theta^{2},
$$
where both $d\sigma^{2}$ and $d\theta^{2}$ are the standard metric of the sphere $S^{2}$.
Let $T, X, \Sigma_{1}, \Sigma_{2}, \Theta_{1}, \Theta_{2}$ be an orthonormal basis of the tangent space
corresponding to the directions $dt, dx, d\sigma, d\theta$ respectively.

Give some constants $c, \gamma, \alpha$ and $t_{1}$ satisfying
$$
0 < c <{\frac 13}, ~~0 < \gamma <{\frac 14}, ~~\alpha > 1, ~~t_{1} > 1.
$$
Set
$$
K=K(c)=\frac{1-c^{2}}{c^{2}},
$$
and define $\psi=\psi(c)$ by
$$
\sin(\sqrt{K}\psi)=\sqrt{1-c^{2}}.
$$
Take $r>0$ satisfying
$$
r \le r(c) = \frac{\pi}{4\sqrt{K}}-\frac{\psi}{2}.
$$
Note that $c$ first is fixed while $r$, $\gamma, \alpha$ and $t_{1}$ will be determined later. We define
$$
t_{i}=t_{1} \alpha^{i},
$$
$$
r_{i}=r t_{i},
$$
$$
\psi_{i}=\psi t_{i},
$$
$$
K_{i}=\frac{K}{t_{i}^{2}},
$$
and
$$
\Delta=\sqrt{K_{i}}(2r_{i}+\psi)=\sqrt{K}(2r+\psi).
$$

In the following, we'll give the constructions of $u, g$ and $f$ respectively.

\subsubsection{Construction of $u(t)$}
For $t_{i} < t < t_{i}+2r_{i}$, we set
$$
u(t)=\frac{1}{\sqrt{K_{i}}} \sin(\sqrt{K_{i}}(t-t_{i}+\psi_{i})),
$$

For $t_{i}+2r_{i} < t < t_{i+1}$, we set
$$
u(t)=t_{i}w(\frac{t}{t_{i}}),
$$
where $w$ is a $C^{2}$ function on $[1+2r,\alpha]$, which is independent of $i$.
Actually, we can define $w$ as follows.  First, set
$$
w(t)=\min\{\frac{\sin\Delta}{\sqrt{K}}+(t-1-2r)\cos\Delta+
\frac{c+1}{2\log\frac{\alpha}{1+2r}}(t\log\frac{t}{1+2r}-t+1+2r), ct \}
$$
and then smoothen $w(t)$ to be a $C^{2}$ function.

Since, at $t=1+2r$,
$$
\frac{\sin\Delta}{\sqrt{K}} < c(1+2r),
$$
and, at $t=\alpha$.
$$
\frac{\sin\Delta}{\sqrt{K}}+(\alpha-1-2r)\cos\Delta+
\frac{c+1}{2\log\frac{\alpha}{1+2r}}(\alpha\log\frac{\alpha}{1+2r}-\alpha+1+2r) > c\alpha,
$$
provided with $\cos\Delta \ge \frac{c+1}{\log\frac{\alpha}{1+2r}}$, which is possible for $\alpha \ge \alpha_{0}(c,r)$. 
(Here, $\alpha_0$ is some constant depending on $c, r$; similarly for the following $\alpha_1, \alpha_2$.)
Thus, $u(t)$ is $C^{1}$ at the endpoints $t_{i}+2r_{i}=(1+2r)t_{i}$ and $t_{i+1}=\alpha t_{i}$.

On the other hand, if $w(t)=ct$,
$$
\begin{cases}
w_{t} \equiv c,\\
w_{tt} \equiv 0;
\end{cases}
$$
if $w(t)=\frac{\sin\Delta}{\sqrt{K}}+(t-1-2r)\cos\Delta+
\frac{c+1}{2\log\frac{\alpha}{1+2r}}(t\log\frac{t}{1+2r}-t+1+2r)$,
$$
\begin{cases}
\cos\Delta \le w_{t} \le \frac{1+3c}{2} < 1,\\
-3\frac{w_{tt}}{w}+\frac{\gamma(1-2\gamma)}{t^{2}} > 0,
\end{cases}
$$
when $\cos(\Delta) \ge \frac{c+1}{\log\frac{\alpha}{1+2r}}$ and $\alpha \ge \alpha_{1}(c,r,\gamma)$.

\vskip .2cm

\noindent
{\bf Conclusion:} When $\alpha \ge \alpha_{1}(c,r,\gamma)$, we have $u(t)$ satisfying, for $t > t_{1}$,
$$
\begin{cases}
\cos\Delta \le u_{t} \le  c, \quad -\frac{u_{tt}}{u}=\frac{K}{t_{i}^{2}}, & t_{i} < t < t_{i}+2r_{i},\\
\cos\Delta \le u_{t} \le \frac{1+3c}{2}, \quad -3\frac{u_{tt}}{u}+\frac{\gamma(1-2\gamma)}{t^{2}} > 0, & t_{i}+2r_{i} < t < t_{i+1}.
\end{cases}
$$

\subsubsection{Construction of $g(t)$}
Now, we construct the function $g(t)$, which is used to control the curvature component
$-\frac{u_{tt}}{u}$ for $t_{i}+2r_{i} < t < t_{i+1}$. More precisely, we let $\gamma<\frac 14$
and $-\frac{g_{tt}}{g} \ge \frac{\gamma(1-2\gamma)}{2t^{2}}$ on these intervals.
On the other hand, in order to make the surgery simple (i.e. we only need to consider surgeries on
$[t_{i}+\frac{r_{i}}{6}, t_{i}+\frac{11r_{i}}{6}] \times_{u(t)} S^3$), we set $g(t)$ to be constants
on $t_{i}+\frac{r_{i}}{6} < t < t_{i}+\frac{11r_{i}}{6}$. Moreover, in order to make
the Ricci curvature $Ric(\Theta_{k},\Theta_{k})$ positive, we need the component
$K(\Theta_{1},\Theta_{2})=\frac{1}{g^{2}}-\frac{g_{t}^{2}}{g^{2}}$.
This is why the dimension of our examples can not be smaller than 6.

For $t_{1} < t < t_{1}+\frac{r_{1}}{6}$, we set
$$
g(t) \equiv g_1 = g(t_{1}+\frac{r_{1}}{6}) = (t_{1}+\frac{r_{1}}{6})^{\gamma} = (1+\frac{r}{6})^{\gamma}t_{1}^{\gamma}.
$$

For $t_{i}+\frac{r_{i}}{6} < t < t_{i+1}+\frac{r_{i+1}}{6}=\alpha(t_{i}+\frac{r_{i}}{6})$, we set
$$
\frac{g_{t}}{g}=\begin{cases}
0, & t_{i}+\frac{r_{i}}{6} < t < t_{i}+\frac{11r_{i}}{6},\\
\frac{6\gamma\beta}{(1+2r)rt_{i}^{2}}[t-(t_{i}+\frac{11r_{i}}{6})], & t_{i}+\frac{11r_{i}}{6} < t < t_{i}+2r_{i},\\
\frac{\gamma\beta}{t}, & t_{i}+2r_{i} < t < t_{i+1},\\
\frac{6\gamma\beta}{\alpha^{2}rt_{i}^{2}}[(t_{i+1}+\frac{r_{i+1}}{6})-t], & t_{i+1} < t < t_{i+1}+\frac{r_{i+1}}{6},
\end{cases}
$$
where $\beta=\frac{\log\alpha}{\log\alpha-\log(1+2r)+\frac{r(1+r)}{6(1+2r)}}=\beta(r,\alpha)$ ensures that
$$
g(t_{i}+\frac{r_{i}}{6})=(t_{i}+\frac{r_{i}}{6})^{\gamma}, \quad \forall i \ge 1.
$$
Thus, for $t_{i}+\frac{r_{i}}{6} < t < t_{i+1}+\frac{r_{i+1}}{6}=\alpha(t_{i}+\frac{r_{i}}{6})$,
$$
g(t) \le g(t_{i+1}+\frac{r_{i+1}}{6})=\alpha^{\gamma}(t_{i}+\frac{r_{i}}{6})^{\gamma} \le \alpha^{\gamma}t^{\gamma}
$$
and
$$
g(t) \ge g(t_{i}+\frac{r_{i}}{6})=(t_{i}+\frac{r_{i}}{6})^{\gamma} \ge \alpha^{-\gamma}t^{\gamma}.
$$
Note that for $\alpha \ge \alpha_{2}(r)$, we can have $\frac 12 \le \beta \le 2$. Then
$$
|\frac{g_{t}}{g}| \le  \frac{2(1+\frac{r}{6})\gamma}{t}
$$
and
$$
\begin{cases}
|\frac{g_{tt}}{g}| \le \frac{12\gamma(1+\frac{r\gamma}{3})}{rt_{i}^{2}}, & t_{i} < t < t_{i}+2r_{i},\\
-\frac{g_{tt}}{g} \ge \frac{\gamma(1-2\gamma)}{2t^{2}}, & t_{i}+2r_{i} < t < t_{i+1}.
\end{cases}
$$

\vskip .2cm

\noindent
{\bf Conclusion:} When $\alpha \ge \alpha_{2}(r)$, we have a $C^{1}$ function $g(t)$ satisfying,  for $t > t_{1}$,
$$
\alpha^{-\gamma}t^{\gamma} \le g(t) \le \alpha^{\gamma}t^{\gamma},
$$
$$
|\frac{g_{t}}{g}| \le  \frac{2(1+\frac{r}{6})\gamma}{t},
$$
and
$$
\begin{cases}
|\frac{g_{tt}}{g}| \le \frac{12\gamma(1+\frac{r\gamma}{3})}{rt_{i}^{2}}, & t_{i} < t < t_{i}+2r_{i},\\
-\frac{g_{tt}}{g} \ge \frac{\gamma(1-2\gamma)}{2t^{2}}, & t_{i}+2r_{i} < t < t_{i+1}.
\end{cases}
$$

\subsubsection{Construction of $f(t,x)$}
The function $f(t,x)$ here is actually $f_{\epsilon(t)}(x)$, which depends on the parameter
$\epsilon<1$. Then, we will set $f_{\epsilon}(x)$ to be almost $R_{0}\sin x$ ($R_{0}$ is sufficiently
small and determined when gluing in $\mathbb{CP}^{2}$'s) and smoothed near $x=0$.
The original idea here comes from \cite{P}, also see \cite{M1} and \cite{M2}.

Let $\phi \in C^{2}(\mathbb{R})$, and satisfy
$$
\phi(x)=\begin{cases}
1, & x \le 0,\\
0, & x \ge 1.
\end{cases}
$$
Set
$$
\phi_\epsilon(x)=\phi\big(\frac{x-\epsilon}{\epsilon^{\frac 14}-\epsilon}\big)=\begin{cases}
1, & x \le \epsilon,\\
0, & x \ge \epsilon^{\frac 14}.
\end{cases}
$$
Let $b=b(\epsilon), l=l(\epsilon)$ and $\delta=\delta(\epsilon)$ to be determined later. Then, set
$$
f_{\epsilon}(x)=\begin{cases}
\frac{\sin(lx)}{l}, & 0 \le x \le b,\\
\frac{\sin(lb)}{l}e^{(1-\epsilon)\log\frac xb}, & b \le x \le \epsilon,\\
R_{0}\sin(x+\delta\phi_{\epsilon}(x)), & \epsilon \le x \le \frac \pi2
\end{cases}
$$
and $f_{\epsilon}(x)=f_{\epsilon}(\pi-x)$ for $\frac \pi2 \le x \le \pi$.

We here remark that $f_\epsilon$ on $[{\frac \pi 2}, \pi]$ is obtained by symmetric extension at $\frac \pi 2$; so,
all the things done around $o_{i}=(t_{i}+r_{i},0)$ in the following
(including the surgery of gluing in $\mathbb{CP}^{2}$'s) can and must also be done
around $o'_{i}=(t_{i}+r_{i},\pi)$ (we have to remove $B_{\frac 45 r_{i}}(o'_{i})$
since the Ricci curvatures on these geodesic balls are not necessarily positive).

Thus, we just consider the case of $0 \le x \le \frac \pi2$ from now on
only by requiring the geodesic balls $B_{\frac 45 r_{i}}(o_{i})$ and $B_{\frac 45 r_{i}}(o'_{i})$ disjoint.
In order to make these two geodesic balls disjoint, we only need to assume
$\frac{\frac 45 r_{i}}{u(t_i)} < \frac{\pi}{2}$, which is valid since we have set
$r \le r(c) = \frac{\pi}{4\sqrt{K}}-\frac{\psi}{2} \le \frac{\pi}{4\sqrt{K}} = \frac{\pi}{4} \frac{c}{\sqrt{1-c^2}}$ and $0<c<\frac 13$.

For the function $f_\epsilon$ in $x$, we have the following two properties,
the proofs of which can be found in \cite{M2}.

\begin{proposition} (\cite{M2}, Lemma 1.22)
$\forall \eta > 0$, there exists $\epsilon_{0} > 0$ such that $\forall \epsilon < \epsilon_{0}, \forall x \in [0,\pi]$,
one has
$$
-\frac{(f_{\epsilon})_{xx}}{f_{\epsilon}} \ge 1-\eta
$$
and
$$
\frac{1-(f_{\epsilon})_{x}^{2}}{f_{\epsilon}^{2}} \ge 1-\eta.
$$
\end{proposition}

\begin{proposition} (\cite{M2}, Lemma 1.26)
There exists $\epsilon_{0} > 0$ such that $\forall \epsilon < \epsilon_{0}$, $\forall x \in [0,\frac \pi2]$,
one has
$$
A(x)=\tan(x)|\frac{(f_{\epsilon})_{x}}{f_{\epsilon}}-\cot(x)| < 2\epsilon.
$$
\end{proposition}

Now, let $\epsilon=\epsilon(t) < \epsilon_{0}$ be a smooth function satisfying
$$
\epsilon(t)=\begin{cases}
\epsilon_{i} & t_{i-1}+\frac 32 r_{i-1} < t < t_{i}+\frac 12 r_{i}\\
\epsilon_{i+1} & t_{i}+\frac 32 r_{i} < t < t_{i+1}+\frac 12 r_{i+1}.
\end{cases}
$$
where $\{\epsilon_i\}$ is a decreasing sequence of positive constants to be determined.

Set $f(t,x)=f_{\epsilon(t)}(x)$. Clearly, $f(t,x)=R_{0}\sin x$ for $\epsilon^{\frac 14} \le x \le \frac \pi2$;
also $f(t,x)=f_{\epsilon_{i}}(x)$ for $t_{i-1}+\frac 32 r_{i-1} < t < t_{i}+\frac 12 r_{i}$ 
and $f(t,x)=f_{\epsilon_{i+1}}(x)$ for $t_{i}+\frac 32 r_{i} < t < t_{i+1}+\frac 12 r_{i+1}$.
This shows that $f_t \neq 0$ is possible only on
$\{(t,x)~|~t_{i}+\frac 12 r_{i} \le t \le t_{i}+\frac 32 r_{i}, 0 \le x \le \epsilon^{\frac 14}\}$.

Thus, if we can show
$\{(t,x)~|~t_{i}+\frac 12 r_{i} \le t \le t_{i}+\frac 32 r_{i},
0 \le x \le \epsilon^{\frac 14}\} \subset B_{\frac 45 r_{i}}(o_{i})$, then $f_{t} \equiv 0$
outside $B_{\frac 45 r_{i}}(o_{i})$.
To prove this, we only need to show
$$
d((t_{i}+r_{i}\pm\frac 12 r_{i},\epsilon^{\frac 14}),(t_{i}+r_{i},0)) < \frac 45 r_{i}.
$$
In fact,
\begin{align*}
    & d((t_{i}+r_{i}\pm\frac 12 r_{i},\epsilon^{\frac 14}),(t_{i}+r_{i},0))\\
\le & d((t_{i}+r_{i}\pm\frac 12 r_{i},\epsilon^{\frac 14}),(t_{i}+r_{i}\pm\frac 12 r_{i},0))
+d((t_{i}+r_{i}\pm\frac 12 r_{i},0),(t_{i}+r_{i},0))\\
\le & \int_{0}^{\epsilon^{\frac 14}} u(t_{i}+\frac 32 r_{i})dx+\frac 12 r_{i}\\
\le & \epsilon^{\frac 14} \frac{t_{i}}{\sqrt{K}}+\frac 12 r_{i}
\end{align*}
Thus, in order of the above aim, we need to choose only $\epsilon < (\frac{r\sqrt{K}}{5})^{4}$.

\vskip .2cm

\noindent
{\bf Conclusion:} For $\epsilon < \epsilon_{0}(c,r)$, $f$ as a function in $x$ always
satisfies the proposition 2.1 and 2.2 above, and $f_{t} \equiv 0$ outside $B_{\frac 45 r_{i}}(o_{i})$.

\vskip .2cm
\noindent {\bf Remark:} {\it We here remark that, at some discrete points, $u(t)$, $g(t)$
and $f_\epsilon (x)$ are only $C^{1}$. However, if the manifold constructed in this manner
has positive Ricci curvature on the complement of those $C^{1}$ parts, the manifold can then
be smoothen to be a $C^{2}$ manifold of positive Ricci curvature; for this, one can refer to \cite{P}, also \cite{M1}.}

\subsection{Construction of $P_i$ ($1 \leq i < +\infty$)}
As mentioned before, each $P_i$ topologically consists of a core and a neck by glueing them together
along the corresponding boundary (component) $S^3$. In this subsection, for sake of convenience,
we introduce a "standard" $P_i'$ (topologically $P_i$) with an appropriate metric according to Perelman so that
we can do surgery on the above doubly warped product; for the details of the proofs, see \cite{P}.
Here, we'd like to remark that by "standard" we do not mean that the metrics have
an explicit or canonical expression, instead it just means that the geometry of the manifolds in question
together with the boundaries has some uniform characteristics (or properties). When we do the surgery, $P_i$ is obtained
by some appropriate rescaling of a certain standard $P_i'$.

The core topologically is $\mathbb{CP}^{2} \setminus B(p_0)$ ($\mathbb{CP}^{2}$ deleting a ball),
but equipped with a metric of positive Ricci curvature
and strictly convex boundary which is isometric to a round sphere. Precisely, this is stated as follows.
\vskip .2cm
\noindent {\bf Core:} {\it The metric can be expressed as
$$
ds^2 = dt^2 + (\sin t \cos t)^2 dx^2 + (\frac {1}{100} \cosh (\frac {t}{100}))^2 dy^2 + (\frac {1}{100} \cosh (\frac {t}{100}))^2 dz^2,
$$
where $0 \leq t \leq t_0$ is the distance from $\mathbb{CP}^{1}$ and $X$,$Y$,$Z$
is the standard invariant framing of $S^3$, satisfying $[X,Y]=2Z$, $[Y,Z]=2X$, $[Z,X]=2Y$.
Moreover, $0<t_{0}<\frac {1}{10}$ satisfies $\frac{\sin(2t_{0})}{2}=\frac{1}{100}\cosh\frac{t_{0}}{100}$,
which makes the boundary a round sphere of radius $\frac{1}{100}\cosh\frac{t_{0}}{100}$
with the normal curvatures bigger than or equal to $\frac{1}{100}\tanh\frac{t_{0}}{100}$.}
\vskip .2cm

The neck is $S^3 \times [0,1]$, equipped with a metric of positive Ricci curvature
such that the boundary component $S^{3} \times \{0\}$ is concave and also isometric to a round sphere,
while the boundary component $S^{3} \times \{1\}$ is convex and looks like a rugby. More precisely, one has the following
\vskip .2cm
\noindent {\bf Neck:} {\it Let $(S^3,g=dt^{2}+B^{2}(t)d\sigma^{2})$
be a rotationally symmetric metric satisfying

{\rm(\romannumeral1)} sectional curvature $> 1$;

{\rm(\romannumeral2)} $0 \le t \le \pi R$, $\max_t\{B(t)\}=r$,
there exists $\rho < \frac{\sinh\frac{t_{0}}{100}}{10^4}$ ($0<t_{0}<\frac {1}{10}$ s.t.
$\frac{\sin(2t_{0})}{2}=\frac{1}{100}\cosh\frac{t_{0}}{100}$) such that $0<r<\rho<R$ and $r^2<\rho^3$.

Then, there exists a metric on $S^3 \times [0,1]$ such that

1) $Ric>0$;

2) the boundary component $S^{3} \times \{0\}$ is concave, with normal curvatures
equal to $-\lambda$, and is isometric to the round sphere $S^{3}(\rho\lambda^{-1})$, for some $\lambda>0$;

3) the boundary component $S^{3} \times \{1\}$ is strictly convex, with all its normal curvatures
bigger than 1, and is isometric to $(S^{3},g)$.}
\vskip .2cm

In order to make the glued manifold have positive Ricci curvature, we need the following criterion, also see \cite{P}.

\vskip .2cm
\noindent {\bf Gluing Criterion:} {\it Let $M_{1}$, $M_{2}$ be two compact smooth manifolds
of positive Ricci curvature, with isometric boundaries $\partial M_{1} \backsimeq \partial M_{2} = X$.
Suppose that the normal curvatures of $\partial M_{1}$ are bigger than the negative of
the normal curvatures of $\partial M_{2}$. Then, the glued manifold $M_{1}\cup_{X} M_{2}$
along the boundary $X$ can be smoothed near $X$ to produce a manifold of positive Ricci curvature.}
\vskip .2cm

According to the above glueing criterion, it is easy to see that, after rescaling the core by $(\frac{1}{100}\cosh\frac{t_{0}}{100})^{-2}\cdot(\frac\rho\lambda)^2$,
we can glue the scaled core (the boundary of which is the round 3-sphere of radius $\rho\lambda^{-1}$ and convex
with the normal curvatures bigger than or equal to $(\frac{\sinh\frac{t_{0}}{100}}{10^4})\cdot\frac\lambda\rho$,
i.e. strictly bigger than $\lambda$) and the neck along the round sphere boundaries to get a "standard" $P_i'$
so that it has positive Ricci curvature and its boundary is isometric to
$(S^3, g)$ and strictly convex with the normal curvatures bigger than $1$. Note that here we only ask $(S^3, g)$
to have the characteristics (or properties) in the definition of the neck.

\subsection{Glue $P_i \times_{g_i} S^2$ with $Q \setminus \coprod\limits_{i=1}^{+\infty}
(B_{\frac 45 r_{i}}(o_{i}) \times_{g_i} S^2)$ ($1 \leq i < +\infty$)}

In \S 2.1, we have constructed $Q=[t_1,+\infty) \times_{u(t)} S^3 \times_{g(t)} S^2$ with the metric
$$
ds^{2}=dt^{2}+u^{2}(t)[dx^{2}+f^{2}(t,x)d\sigma^{2}]+g^{2}(t)d\theta^{2}.
$$
Moreover, we also mentioned that we only need to do surgeries
on $(t_{i}+\frac{r_{i}}{6}, t_{i}+\frac{11r_{i}}{6}) \times_{u(t)} S^3$,
since $g(t)\equiv g_i$ (a constant) on $[t_{i}+\frac{r_{i}}{6}, t_{i}+\frac{11r_{i}}{6}]$.

In order to do surgery, we first need to remove a geodesic ball $B_{\frac 45 r_{i}}(o_{i})$
from $(t_{i}+\frac{r_{i}}{6}, t_{i}+\frac{11r_{i}}{6}) \times_{u(t)} S^3$.
Note that $u(t)=\frac{1}{\sqrt{K_{i}}} \sin(\sqrt{K_{i}}(t-t_{i}+\psi_{i}))$ on $[t_{i}, t_{i}+2r_i]$.

Then, we need to check that the boundary $\partial B_{\frac 45 r_{i}}(o_{i})$
with the induced metric, after some suitable rescaling, owns the metric properties of the boundary
of a certain standard $P_i'$, i.e. the assumptions (i) and (ii) of the metric $g$ in the definition
of the neck in \S 2.2, so that we can glue the $P_i$ (obtained from such a $P_i'$ by a suitable scaling)
on $((t_{i}+\frac{r_{i}}{6}, t_{i}+\frac{11r_{i}}{6}) \times_{u(t)} S^3) \setminus B_{\frac 45 r_{i}}(o_{i})$
along the isometric boundaries $\partial B_{\frac 45 r_{i}}(o_{i})$ and $\partial P_i$.

Furthermore, according to the gluing criterion in \S 2.2, in order to guarantee the resulted manifold
on $(t_{i}+\frac{r_{i}}{6}, t_{i}+\frac{11r_{i}}{6})$
to admit a metric with positive Ricci curvature, we need to check that the Ricci curvature of
each piece of manifold is positive and the normal curvatures of $\partial P_i$ are bigger
than the negatives of the normal curvatures of $\partial B_{\frac 45 r_{i}}(o_{i})$
in $((t_{i}+\frac{r_{i}}{6}, t_{i}+\frac{11r_{i}}{6}) \times_{u(t)} S^3) \setminus B_{\frac 45 r_{i}}(o_{i})$.

Obviously, when we take the product of the above resulted manifold with a two-sphere (i.e. $\times_{g_i} S^2$, $g_i$ is constant),
it still has positive Ricci curvature. Thus, if both
$Q \setminus \coprod\limits_{i=1}^{+\infty} (B_{\frac 45 r_{i}}(o_{i}) \times_{g_i} S^2)$ and
$P_i \times_{g_i} S^2$ have positive Ricci curvature (actually, the later is already known of positive Ricci curvature in \S 2.2),
they can be glued together to get a manifold with positive Ricci curvature. The positiveness of Ricci curvature of
$Q \setminus \coprod\limits_{i=1}^{+\infty} (B_{\frac 45 r_{i}}(o_{i}) \times_{g_i} S^2)$ will be verified in \S 3.

In the remaining part of this subsection, we shall check the following two facts.
The first one is that the normal curvatures of $\partial P_i$ are bigger than
the negatives of the normal curvatures of $\partial B_{\frac 45 r_{i}}(o_{i})$.
The second one is that the boundary $\partial B_{\frac 45 r_{i}}(o_{i})$ with the induced metric,
after some suitable rescaling, owns the characteristics of the metric $g$ in the definition of the neck in \S 2.2.

\subsubsection{Normal curvatures}

In order to prove the first fact, it is sufficient to show that for $\epsilon$ sufficiently small,
any $x \in \partial B_{\frac 45 r_{i}}(o_{i})$, and any orthogonal unit vectors $X,Y\in S_x \partial B_{\frac 45 r_{i}}(o_{i})$,
$$
K_{int}(X,Y) > \max_{Z \in S_x \partial B_{\frac 45 r_{i}}(o_{i})} |\uppercase\expandafter{\romannumeral2}_{N}(Z,Z)|^2,
$$
where $N$ is the unit outward normal vector of $\partial B_{\frac 45 r_{i}}(o_{i})$,
$\uppercase\expandafter{\romannumeral2}_{N}$ is the corresponding second fundamental form,
$K_{int}$ is the intrinsic curvature, and
$S\partial B_{\frac 45 r_{i}}(o_{i}) \subset T\partial B_{\frac 45 r_{i}}(o_{i})$ is the unit tangent bundle.

In fact, after rescaling $\partial B_{\frac 45 r_{i}}(o_{i})$ by
$\max_{Z \in S_x \partial B_{\frac 45 r_{i}}(o_{i})} |\uppercase\expandafter{\romannumeral2}_{N}(Z,Z)|)^2$,
we obtain a rotationally symmetric metric $(S^{3},g)$ with curvature bigger than 1 and
the absolute value of the corresponding normal curvatures less than or equal to 1. In the next subsection,
we will check that this rotationally symmetric metric owns the characteristics of the metric $g$ of the neck in \S 2.2,
so that it is the boundary of a certain standard $P_i'$.
In particular, the normal curvatures of $\partial P_i$ in $P_i$ (obtained from the $P'_i$ by
$(\max_{Z \in S_x \partial B_{\frac 45 r_{i}}(o_{i})} |\uppercase\expandafter{\romannumeral2}_{N}(Z,Z)|)^{-2}$-rescaling)
are bigger than the negative of the normal curvatures of $\partial B_{\frac 45 r_{i}}(o_{i})$.

The proof of the above inquality is similar to that in \S1.6 of \cite{M1}.
Denote by $N$ the unit outward normal vector to $\partial B_{\frac 45 r_{i}}(o_{i})$,
$$
N=T\cos\xi+X\sin\xi,
$$
and
$$
Y=X\cos\xi-T\sin\xi.
$$
Thus
$$
\uppercase\expandafter{\romannumeral2}_{N}(Y,Y)=\sqrt{K_{i}}\cot(\frac 45 \sqrt{K_{i}}r_{i}),
$$
since the $T\wedge X$ plane is isometric to the corresponding one of $S^{2}(\frac{1}{\sqrt{K_{i}}})$.
On the other hand, for $j=1, 2$,
$$
\uppercase\expandafter{\romannumeral2}_{N}(\Sigma_{j},\Sigma_{j})=\frac{u_{t}}{u}\cos\xi+\frac{f_{x}}{fu}\sin\xi.
$$

Note that if $f(t,x)=\sin x$,
$$
\frac{u_{t}}{u}\cos\xi+\frac{f_{x}}{fu}\sin\xi=\sqrt{K_{i}}\cot(\frac 45 \sqrt{K_{i}}r_{i}).
$$
Thus,
\begin{align*}
|\uppercase\expandafter{\romannumeral2}_{N}(\Sigma_{j},\Sigma_{j})-\sqrt{K_{i}}\cot(\frac 45 \sqrt{K_{i}}r_{i})| & \le A(x)\frac{\cot x |\sin\xi|}{u}\\
                    & = A(x)|\sqrt{K_{i}}\cot(\frac 45 \sqrt{K_{i}r_{i}})-\frac{u_{t}}{u}\cos\xi|\\
                    & \le A(x)(\sqrt{K_{i}}\cot(\frac 45 \sqrt{K_{i}r_{i}})+|\frac{u_{t}}{u}|)\\
                    & \le 2\epsilon[\frac{\sqrt{K}\cot(\frac 45 \sqrt{K}r)}{t_{i}}+\frac{1+3c}{2ct_{i}}]\\
                    & = 2[\sqrt{K}\cot(\frac 45 \sqrt{K}r)+\frac{1+3c}{2c}]\frac{\epsilon}{t_{i}}\\
                    & = D(c,r)\frac{\epsilon}{t_{i}} \le D(c,r)\frac{\epsilon_i}{t_i}.
\end{align*}

Then, one has
$$
\frac{\sqrt{K}\cot(\frac 45 \sqrt{K}r)}{t_{i}}\leq\max_{Z \in S_x \partial B_{\frac 45 r_{i}}(o_{i})} |\uppercase\expandafter{\romannumeral2}_{N}(Z,Z)|\leq\frac{\sqrt{K}\cot(\frac 45 \sqrt{K}r)}{t_{i}}+D(c,r)\frac{\epsilon_i}{t_i}.
$$

For the intrinsic curvatures, by Gauss equation, we have
\begin{align*}
K_{int}(\Sigma_{1},\Sigma_{2}) & = K(\Sigma_{1},\Sigma_{2})+(\uppercase\expandafter{\romannumeral2}_{N}(\Sigma_{j},\Sigma_{j}))^{2}\\
                             & \ge \frac{1-f_{x}^{2}}{f^{2}}\frac{1}{u^{2}}-\frac{u_{t}^{2}}{u^{2}}+[\frac{\sqrt{K}\cot(\frac 45 \sqrt{K}r)}{t_{i}}-D(c,r)\frac{\epsilon}{t_{i}}]^{2}\\
                             & \ge \frac{1-\eta-(\frac{1+3c}{2})^{2}}{u^{2}}+[\frac{\sqrt{K}\cot(\frac 45 \sqrt{K}r)}{t_{i}}-D(c,r)\frac{\epsilon}{t_{i}}]^{2}\\
                             & \ge \frac{1-\eta-(\frac{1+3c}{2})^{2}}{c^{2}t_{i}^{2}}+[\frac{\sqrt{K}\cot(\frac 45 \sqrt{K}r)}{t_{i}}-D(c,r)\frac{\epsilon}{t_{i}}]^{2}\\
                             & > [\frac{\sqrt{K}\cot(\frac 45 \sqrt{K}r)}{t_{i}}+D(c,r)\frac{\epsilon_i}{t_i}]^{2},
\end{align*}
assuming $\epsilon$ sufficiently small. Similarly,
\begin{align*}
K_{int}(Y,\Sigma_{j}) & = K(X,\Sigma_{j})\cos^{2}\xi+K(T,\Sigma_{j})\sin^{2}\xi+(\uppercase\expandafter{\romannumeral2}_{N}(Y,Y)
(\uppercase\expandafter{\romannumeral2}_{N}(\Sigma_{j},\Sigma_{j}))\\
                             & \ge (-\frac{f_{xx}}{f}\frac{1}{u^{2}}-\frac{u_{t}^{2}}{u^{2}})\cos^{2}\xi+
                             (-\frac{u_{tt}}{u})\sin^{2}\xi\\
                             & \quad +\frac{\sqrt{K}\cot(\frac 45 \sqrt{K}r)}{t_{i}}[\frac{\sqrt{K}\cot(\frac 45 \sqrt{K}r)}{t_{i}}-D(c,r)\frac{\epsilon}{t_{i}}]\\
                             & \ge \frac{1-\eta-(\frac{1+3c}{2})^{2}}{c^{2}t_{i}^{2}}\cos^{2}\xi+
                             \frac{K}{t_{i}^{2}}\sin^{2}\xi\\
                             & \quad +\frac{\sqrt{K}\cot(\frac 45 \sqrt{K}r)}{t_{i}}[\frac{\sqrt{K}\cot(\frac 45 \sqrt{K}r)}{t_{i}}-D(c,r)\frac{\epsilon}{t_{i}}]\\
                             & > [\frac{\sqrt{K}\cot(\frac 45 \sqrt{K}r)}{t_{i}}+D(c,r)\frac{\epsilon_i}{t_i}]^{2},
\end{align*}
assuming $\epsilon$ sufficiently small.

\subsubsection{The boundary $\partial B_{\frac 45 r_{i}}(o_{i})$ with the induced metric}

Next, we prove the second fact, i.e. $\partial B_{\frac 45 r_{i}}(o_{i})$ with the induced metric,
after rescaling by $\max_{Z \in S_x \partial B_{\frac 45 r_{i}}(o_{i})} |\uppercase\expandafter{\romannumeral2}_{N}(Z,Z)|)^2$,
owns the characteristics of the metric $g$ of the neck in \S 2.2.

Even though we could not write out the precise expression of the metric, we can check it in the limit case (thanks to the construction of the neck by Perelman \cite{P}, it only requires some rough characterization of the geometric shape).

Indeed, for $i$ sufficiently large, as $\epsilon$ goes to 0, the metric of $B_{\frac 45 r_{i}}(o_{i})$ converges to
$$
dt^{2}+(\frac{1}{\sqrt{K_{i}}} \sin(\sqrt{K_{i}}(t-t_{i}+\psi_{i})))^{2}[dx^{2}+(R_{0}\sin x)^{2}d\sigma^{2}].
$$
Thus the metric of $\partial B_{\frac 45 r_{i}}(o_{i})$ converges to
$$
(\frac{t_{i}\sin(\frac 45 \sqrt{K}r)}{\sqrt{K}})^{2}[dy^{2}+(R_{0}\sin y)^{2}d\sigma^{2}].
$$
Rescaled by $(\frac{\sqrt{K}\cot(\frac 45 \sqrt{K}r)}{t_{i}})^2$ (note that $\max_{Z \in S_x \partial B_{\frac 45 r_{i}}(o_{i})} |\uppercase\expandafter{\romannumeral2}_{N}(Z,Z)|)^2$ converges to $(\frac{\sqrt{K}\cot(\frac 45 \sqrt{K}r)}{t_{i}})^2$ when $\epsilon$ goes to 0) and letting $t=y\cos(\frac 45 \sqrt{K}r)$, the metric becomes
$$
dt^{2}+(R_{0}\cos(\frac 45 \sqrt{K}r)\sin\frac{t}{\cos(\frac 45 \sqrt{K}r)})^{2}d\sigma^{2}.
$$
Choosing $R_{0}$ sufficiently small, the metric is exactly what we want.

\subsection{Construction of $P_0$}

Note that after gluing $P_i$ ($1 \leq i < +\infty$) with $Q$, the manifold has a boundary,
i.e. the hypersurface $t=t_{1}$. In this subsection, we will glue in a core via a neck along the boundary $t=t_1$
to smooth the manifold near the origin, based on the gluing criterion.

For $t_{1} < t < t_{1}+\frac{r_{1}}{6}$, we already have
$$
\begin{cases}
u(t) = \frac{1}{\sqrt{K_{1}}} \sin(\sqrt{K_{1}}(t-t_{1}+\psi_{1})),\\
g(t) \equiv (t_{1}+\frac{r_{1}}{6})^{\gamma}.
\end{cases}
$$
Then, for the hypersurface $t=t_{1}$, the intrinsic metric is
$$
ds^{2}=(ct_{1})^{2}[dx^{2}+f_{\epsilon_{1}}^{2}(x)d\sigma^{2}].
$$
Rescaled by $t_{1}^{-2}$ and letting $t=cx$, when $\epsilon_{1}$ goes to 0, it converges to
$$
dt^{2}+(cR_{0}\sin\frac{t}{c})^{2}d\sigma^{2}.
$$
Similarly, choosing $R_{0}$ sufficiently small, the metric owns the properties of the metric $g$ of the neck in \S 2.2.

On the other hand,
$$
\uppercase\expandafter{\romannumeral2}_{-T}(X,X)=
\uppercase\expandafter{\romannumeral2}_{-T}(\Sigma_{j},\Sigma_{j})=-\frac{1}{t_{1}},
$$
$$
K_{int}(\Sigma_{1},\Sigma_{2})=K(\Sigma_{1},\Sigma_{2})+\frac{1}{t_{1}^{2}}=
\frac{1-f_{x}^{2}}{f^{2}}\frac{1}{u^{2}} \ge \frac{1-\eta}{c^{2}t_{1}^{2}},
$$
$$
K_{int}(X,\Sigma_{j})=K(X,\Sigma_{j})+\frac{1}{t_{1}^{2}}=
-\frac{f_{xx}}{f}\frac{1}{u^{2}} \ge \frac{1-\eta}{c^{2}t_{1}^{2}}.
$$
Setting $1-\eta > c^2$, we can continue the manifold for $t \le t_{1}$ by gluing a core via a neck (i.e. $P_0$) with a metric product $\times_{g_0}S^{2}$ ($g_0=(t_{1}+\frac{r_{1}}{6})^{\gamma}$).

\section{Positiveness of Ricci curvature and quadratically asymptotic non-negativeness of curvature}

\subsection{Positive Ricci curvature}

In this subsection, we will show that the constructed manifold above has positive Ricci curvature
by taking some appropriate constants, i.e. $\eta, r, \gamma$ sufficiently small and $\alpha, t_1$ sufficiently big.

Recall that our manifold is constructed as
$$
M = (Q \setminus \coprod\limits_{i=1}^{+\infty} (B_{\frac 45 r_{i}}(o_{i}) \times_{g_i} S^2))
\cup_{\textrm{Id}} \coprod\limits_{i=0}^{+\infty} (P_i \times_{g_i} S^2).
$$
As mentioned before, by our constrction and the gluing criterion, we only need to check that
$Q \setminus \coprod\limits_{i=1}^{+\infty} (B_{\frac 45 r_{i}}(o_{i}) \times_{g_i} S^2)$ have positive Ricci curvature.

Since $f_{t} \equiv 0$ outside $B_{\frac 45 r_{i}}(o_{i})$, the nonzero Ricci curvatures are as follows.
\begin{align*}
Ric(T,T) & = -3\frac{u_{tt}}{u}-2\frac{g_{tt}}{g}\\
         & \ge \begin{cases}
         -3\frac{u_{tt}}{u}+\frac{\gamma(1-2\gamma)}{t^{2}} > 0, & t_{i}+2r_{i} < t < t_{i+1},\\
         3\frac{K}{t_{i}^{2}}-\frac{24\gamma(1+\frac{r\gamma}{3})}{rt_{i}^{2}} > 0, & t_{i} < t < t_{i}+2r_{i},
         \end{cases}
\end{align*}
provided with $\gamma < \gamma_{0}(c,r)$;

\begin{align*}
Ric(X,X) & =-2\frac{f_{xx}}{f}\frac{1}{u^{2}}-2(\frac{u_{t}}{u})^{2}-2\frac{u_{t}}{u}\frac{g_{t}}{g}-\frac{u_{tt}}{u}\\
         & \ge 2\frac{1-\eta}{u^{2}}-2(\frac{u_{t}}{u})^{2}--2\frac{u_{t}}{u}\frac{2(1+\frac r6)\gamma}{t}-\frac 13 \frac{\gamma(1-2\gamma)}{t^{2}}\\
         & = \frac{2}{u^{2}}[1-\eta-u_{t}^{2}-\frac{2(1+\frac r6)\gamma u_{t}u}{t}]-\frac 13 \frac{\gamma(1-2\gamma)}{t^{2}}\\
         & \ge \frac{2}{(\frac{1+3c}{2})^{2}t^{2}}\{1-\eta-(\frac{1+3c}{2})^{2}[1+2(1+\frac r6)\gamma]\}-\frac 13 \frac{\gamma(1-2\gamma)}{t^{2}}\\
         & = \frac{2}{(\frac{1+3c}{2})^{2}t^{2}}\{1-\eta-(\frac{1+3c}{2})^{2}[1+2(1+\frac r6)\gamma+\frac 16 \gamma(1-2\gamma)]\}\\
         & > 0,
\end{align*}
provided with $1-\eta > (\frac{1+3c}{2})^{2}[1+2(1+\frac r6)\gamma+\frac 16 \gamma(1-2\gamma)]$,
which is possible when $\eta < \eta_{0}(c)$ and $\gamma < \gamma_{0}(c,r,\eta)$;

\begin{align*}
Ric(\Sigma_{j},\Sigma_{j}) & =(\frac{1-f_{x}^{2}}{f^{2}}-\frac{f_{xx}}{f})\frac{1}{u^{2}}-2(\frac{u_{t}}{u})^{2}-
2\frac{u_{t}}{u}\frac{g_{t}}{g}-\frac{u_{tt}}{u},\\
                           & > 0
\end{align*}
which is totally similar to $Ric(X,X)$;

\begin{align*}
Ric(\Theta_{k},\Theta_{k}) & =\frac{1}{g^{2}}-(\frac{g_{t}}{g})^{2}-\frac{g_{tt}}{g}-3\frac{u_{t}}{u}\frac{g_{t}}{g}\\
                           & \ge \frac{1}{\alpha^{\gamma}t^{\gamma}}-\frac{4(1+\frac r6)^{2}\gamma^{2}}{t^{2}}-\frac{12\gamma(1+2r)(1+\frac{r\gamma}{3})}{rt^{2}}-
                           \frac{3(1+3c)(1+\frac r6)\gamma}{t^{2}\cos\Delta}\\
                           & > 0,
\end{align*}
provided with $t_{1} > t_{1}(c,r,\gamma)$.

Moreover, all the off-diagonal terms of the Ricci tensor vanish. Thus the Ricci curvatures are positive.




\subsection{Quadratically asymptotic non-negativeness of sectional curvature}
In this subsection, we show that our example has quadratically asymptotically nonnegative sectional curvature.

For $P_i \times_{g_i} S^2$, we only have to check that $P_i$ has quadratically asymptotically nonnegative sectional curvature. After rescaling the standard cores and necks by $\max_{Z \in S_x \partial B_{\frac 45 r_{i}}(o_{i})} |\uppercase\expandafter{\romannumeral2}_{N}(Z,Z)|)^2$ ($\leq(\frac{\sqrt{K}\cot(\frac 45 \sqrt{K}r)}{t_{i}}+D(c,r)\frac{\epsilon_i}{t_i})^2$),
$$
K_{p_{0}} (t) \ge -(K_{0})^{2}(\frac{\sqrt{K}\cot(\frac 45 \sqrt{K}r)}{t_{i}}+
D(c,r)\frac{\epsilon}{t_{i}})^{2} \ge -\frac{K_{0}^{2}(c,r)}{t_{i}^{2}}.
$$

For $Q \setminus \coprod\limits_{i=1}^{+\infty} (B_{\frac 45 r_{i}}(o_{i}) \times_{g_i} S^2)$,
$$
K(T,X)=K(T,\Sigma_{j})=-\frac{u_{tt}}{u} > -\frac 13 \frac{\gamma(1-2\gamma)}{t^{2}},
$$
$$
K(T,\Theta_{k})=-\frac{g_{tt}}{g} \ge -\frac{12\gamma(1+\frac{r\gamma}{3})(1+2r)}{rt^{2}},
$$
$$
K(X,\Sigma_{j})=-\frac{f_{xx}}{f}\frac{1}{u^{2}}-\frac{u_{t}^{2}}{u^{2}} \ge \frac{1-\eta-(\frac{1+3c}{2})^{2}}{\cos^{2}\Delta}\frac{1}{t^{2}} > 0,
$$
$$
K(X,\Theta_{k})=K(\Sigma_{j},\Theta_{k})=-\frac{u_{t}}{u}\frac{g_{t}}{g} \ge -\frac{(1+3c)(1+\frac r6)\gamma}{\cos\Delta}\frac{1}{t^{2}},
$$
$$
K(\Sigma_{1},\Sigma_{2})=\frac{1-f_{x}^{2}}{f^{2}}\frac{1}{u^{2}}-\frac{u_{t}^{2}}{u^{2}} \ge \frac{1-\eta-(\frac{1+3c}{2})^{2}}{\cos^{2}\Delta}\frac{1}{t^{2}} > 0,
$$
$$
K(\Theta_{1},\Theta_{2})=\frac{1}{g^{2}}-\frac{g_{t}^{2}}{g^{2}} \ge \frac{1}{t^{2\gamma}}-\frac{4(1+\frac r6)^{2}\gamma^{2}}{t^{2}} > 0.
$$
And other terms of curvature tensors are zero.

Thus, we have
$$
K_{p_{0}} (t) \ge -\frac{K_{0}^{2}(c,r,\gamma)}{t^{2}}.
$$

\vskip .2cm
\noindent {\bf Remark:} {\it
For the diameter growth, we have
$$
0 < \pi\cos\Delta = \liminf_{t \to \infty} \frac{{\text{diam}}(p_0; t)}{t}
< \limsup_{t \to \infty} \frac{{\text{diam}}(p_0; t)}{t} = \pi c.
$$
For the volume growth, we have
$$
0 < V_{1} = \liminf_{t \to \infty} \frac{{\text{Vol}}(B(p_0,t))}{t^{4+2\gamma}}
< \limsup_{t \to \infty} \frac{{\text{Vol}}(B(p_0,t))}{t^{4+2\gamma}} = V_{2} < +\infty.
$$}

\vskip.5cm
\noindent School of Mathematical Sciences, Shanghai Jiao Tong University, Shanghai, 200240\\
{\it E-mail}: dan@sjtu.edu.cn
\vskip .3cm
\noindent School of Mathematical Sciences, Shanghai Jiao Tong University, Shanghai, 200240\\
{\it E-mail}: {yangyihu@sjtu.edu.cn}

\end{document}